\documentclass[10pt,letter]{article}

\usepackage{amsmath}
\usepackage{amsfonts}
\usepackage{amssymb}
\usepackage{graphicx}
\usepackage{enumerate}

\begin{document}

\newtheorem{theorem}{Theorem}[section]
\newtheorem{corollary}[theorem]{Corollary}
\newtheorem{definition}[theorem]{Definition}
\newtheorem{conjecture}[theorem]{Conjecture}
\newtheorem{question}[theorem]{Question}
\newtheorem{lemma}[theorem]{Lemma}
\newtheorem{proposition}[theorem]{Proposition}
\newtheorem{example}[theorem]{Example}
\newenvironment{proof}{\noindent {\bf
Proof.}}{\rule{3mm}{3mm}\par\medskip}
\newcommand{\remark}{\medskip\par\noindent {\bf Remark.~~}}
\newcommand{\pp}{{\it p.}}
\newcommand{\de}{\em}

\title{ {\bf The Harary index of trees}
\thanks{
Supported by the Research Grant 144007 of Serbian Ministry of Science, the Postdoctoral Science
Foundation of Central South University, China Postdoctoral Science Foundation, (No. 70901048,
10871205) and NSFSD (No. Y2008A04, BS2010SF017). }}

\author{
Aleksandar Ili\' c \footnotemark[3]\\
{\small Faculty of Sciences and Mathematics, University of Ni\v s} \\
{\small Vi\v segradska 33, 18000 Ni\v s, Serbia} \\
{\small e-mail: { \tt aleksandari@gmail.com}}\vspace{0.3cm} \\
Guihai Yu\\
{\small Department of Mathematics, Shandong Institute of Business and Technology} \\
{\small 191 Binhaizhong Road, Yantai, Shandong, P.R. China, 264005.}\\
{\small e-mail: {\tt yuguihai@126.com }} \vspace{0.3cm} \\
Lihua Feng \\
{\small Department of Mathematics, Central South University} \\
{\small Railway Campus,  Changsha, Hunan, P. R. China, 410075.}\\
{\small e-mail: {\tt fenglh@163.com }} \\
}
\maketitle \vspace{-0.6cm}

\begin{abstract}
The Harary index of a graph $G$ is recently introduced topological index, defined on the reverse
distance matrix as $H(G)=\sum_{u,v \in V(G)}\frac{1}{d(u,v)}$, where $d(u,v)$ is the length of the
shortest path between two distinct vertices $u$ and~$v$. We present the partial ordering of
starlike trees based on the Harary index and we describe the trees with the second maximal and the
second minimal Harary index. In this paper, we investigate the Harary index of trees with $k$
pendent vertices and determine the extremal trees with maximal Harary index. We also characterize
the extremal trees with maximal Harary index with respect to the number of vertices of degree two,
matching number, independence number, radius and diameter. In addition, we characterize the
extremal trees with minimal Harary index and given maximum degree. We concluded that in all
presented classes, the trees with maximal Harary index are exactly those trees with the minimal
Wiener index, and vice versa.
\end{abstract}


 \vskip 0.1cm

\section{Introduction}
\vskip -0.1cm

 In theoretical chemistry molecular structure descriptors (also called topological
indices) are used for modeling physico-chemical, pharmacologic, toxicologic, biological and other
properties of chemical compounds \cite{todeschini, ToCo1}. There exist several types of such
indices, especially those based on graph-theoretical distances.

In 1993 Plav\v {s}i\'c et al. in \cite{plavsic} and Ivanciuc et al. in \cite{ivanciuc1}
independently introduced a new topological index, which was named Harary index in honor of Frank
Harary on the occasion of his 70th birthday. This topological index is derived from the reciprocal
distance matrix and has a number of interesting chemical-physics properties \cite{ivanciuc3}. The
Harary index and its related molecular descriptors have shown some success in structure-property
correlations \cite{devillers,diudea,diudea1,gutman,ivanciuc2,lucic,zhou0}. Its modification has
also been proposed \cite{lucic2} and their use in combination with other molecular descriptors
improves the correlations \cite{todeschini, trinajstic}. It is of interest to study spectra and
polynomials of these matrices \cite{GuKlYaYe06,Zhou4}.

In this paper, let $G$ be a simple connected (molecular) graph with
vertex set $V(G)$. The Harary index is defined as the half-sum of the elements
in the reciprocal distance matrix (also called the Harary matrix \cite{janezic}),
$$
H(G)=\sum_{u,v \in V(G)}\frac{1}{d(u,v)},
$$
where $d(u,v)$ is the  distance between $u$ and $v$ in $G$ and the
sum goes over all the pairs of vertices. The Wiener index, defined
as
$$
W (G) = \sum_{u, v \in V (G)} d (u, v),
$$
is considered as one of the most used topological indices with high correlation with many physical
and chemical properties of molecular compounds. The  majority of chemical applications of the
Wiener index deal with acyclic organic molecules. For recent results and applications of Wiener
index see~\cite{DoEn01}.

Up to now, many results were obtained concerning the Harary index of a graph. Gutman~\cite{gutman}
supported the use of Harary index as a measure of branching in alkanes, by showing

\begin{theorem}
Let $T$ be a tree on $n$ vertices. Then
$$
1+n \sum_{k=2}^{n-1} \frac 1 k \leq H(T) \leq \frac{(n+2)(n-1)}{4}.
$$
 The right  equality holds if and only if
$T\cong S_{n}$, while the left equality holds if and only if $T\cong
P_{n}$.
\end{theorem}

In this paper, we further refine this relation by introducing long chain of inequalities and obtain
the trees with the second maximum and the second minimum Harary index. In \cite{zhou1} Zhou, Cai
and Trinajsti\'{c} presented some lower and upper bounds for the Harary index of connected graphs,
triangle-free and quadrangle-free graphs and gave the Nordhaus-Gaddum-type inequalities. In
\cite{das} the authors obtained some lower and upper bounds for the Harary index of graphs in terms
of the diameter and the number of connected components. Zhou, Du and Trinajsti\' c \cite{zhou2}
discussed the Harary index of landscape graphs, which have found applications in ecology
\cite{UrKe02}. Feng and Ili\' c in \cite{IlFe10} establish sharp upper bounds for the Zagreb
indices, sharp upper bound for the Harary index and sharp lower bound for the hyper-Wiener index of
graphs with a given matching number.

In this paper, we analyze relations between the extremal trees with maximum and minimum Harary and
Wiener index. We investigate the Harary index of $n$-vertex trees with given number of pendent
vertices and determine the extremal trees with maximal $H (G)$. Furthermore, we derive the partial
ordering of starlike trees based on the majorization inequalities of the pendent path lengths. We
characterize the extremal trees with maximal Harary index in terms of the number of vertices of
degree two, the matching number, independence number, radius and diameter. Finally, we characterize
the extremal trees with minimal Harary index with respect to the maximum vertex degree. All these
results are compared with those of the ordinary Wiener index. We conclude the paper by posing a
conjecture regarding to the extremal tree with the maximum Harary index among $n$-vertex trees with
fixed maximum degree.

\section{Preliminaries}
\vskip -0.1cm

For any two vertices $u$ and $v$ in $G$, the distance between $u$ and $v$, denoted by $d_{G}(u,v),$
is the number of edges in a shortest path joining $u$ and $v$. The eccentricity $\varepsilon (v)$
of a vertex $v$ is the maximum distance from $v$ to any other vertex. The vertices of minimum
eccentricity form the center. A tree has exactly one or two adjacent center vertices; in this
latter case one speaks of a bicenter. The diameter $d (G)$ of a graph G is the maximum eccentricity
over all vertices in a graph, while the radius $r (G)$ is the minimum eccentricity over all $v \in
V (G)$. For a vertex $u$ in $G$, the degree of $u$ is denoted by $deg(u)$.

Two distinct edges in a graph $G$ are {\it independent} if they are not incident with a common
vertex in $G$. A set of pairwise independent edges in $G$ is called a {\it matching} in $G$, while
a matching  of maximum cardinality is a {\it maximum matching} in $G$. The {\it matching number}
$\beta(G)$ of $G$ is the cardinality of a maximum matching of $G$. It is well known that $\beta(G)
\leq \frac{n}{2}$, with equality if and only if $G$ has a perfect matching. The {\it independence
number} of $G$, denoted by $\alpha(G)$, is the size of a maximum independent set of $G$.



Let $P_n$ and $S_n$ denote the path and the star on $n$ vertices.
 A {\it starlike tree} is a tree with exactly
one  vertex of degree at least 3. We denote by $S
(n_{1},n_{2},\ldots,n_{k})$ the starlike tree  of order $n$ having a
branching vertex $v$ and
$$
S (n_{1},n_{2},\ldots,n_{k})-v=P_{n_1}\cup P_{n_2}\cup \ldots \cup P_{n_k},
$$
where $n_1\geq n_2\geq \ldots\geq n_k \geq 1$. Clearly, the numbers $n_1, n_2,
\ldots, n_k$ determine the starlike tree up to isomorphism and $n = n_1 + n_2 + \ldots + n_k + 1$.
The starlike tree $BS_{n, k} \cong S (n_{1},n_{2},\ldots,n_{k})$ is {\em balanced}
if all paths have almost equal lengths, i.e., $|n_i - n_j| \leqslant 1$ for every $1
\leqslant i < j \leqslant k$.


Denote by $\Delta (T)$ the maximum vertex degree of a tree $T$. The path $P_n$ is the unique tree
with $\Delta = 2$; while the star $S_n$ is the unique tree with $\Delta = n-1$. Therefore, we can
assume that $3 \leq \Delta \leq n - 2$. The broom $B_{n, \Delta}$ is a tree consisting of a star
$S_{\Delta + 1}$ and a path of length $n - \Delta - 1$ attached to an arbitrary pendent vertex of
the star.

If $\frac{n-1}{2} < m \leq n-1$, then $A_{n,m}$ is the tree obtained from $S_{m+1}$ by adding a
pendent edge to each of $n-m-1$ of the the pendent vertices of $S_{m+1}$. We call $A_{n,m}$ a spur
(see Fig. 1). Clearly, $A_{n,m}$ has $n$ vertices and $m$ pendent vertices; the matching number and
the independence number of $A_{n,m}$ are $n-m$ and $m$, respectively. Note that if $m >
\frac{n-1}{2}$, then $A_{n,m}\cong BS_{n,m}$.

\begin{figure}[h]
  \center
  \includegraphics [width = 4cm]{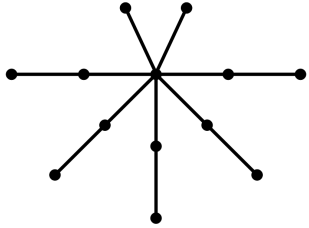}
  \caption { \textit{ The spur $A_{13, 7}$\,. } }
\end{figure}

Next, we give some lemmas which are very useful in the following.

\begin{definition}
\label{de-delta} Let $v$ be a vertex of a tree $T$  and $deg(v) =m + 1$. Suppose that $P_1, P_2,
\ldots, P_m$ are pendent paths incident with $v$, with  the starting vertices of paths $v_1, v_2,
\ldots, v_m$, respectively and lengths $n_i \geqslant 1$ $(i = 1, 2, \ldots, m)$. Let $w$ be the
neighbor of $v$ distinct from  $v_i$. Let $T' = \delta (T, v)$ be a tree obtained from $T$ by
removing the edges $v v_1, v v_2, \ldots, v v_{m - 1}$
 and adding  edges $w v_1, w v_2, \ldots, w v_{m
- 1}$. We say that $T'$ is a $\delta$-transform of~$T$.
\end{definition}

This transformation preserves the number of pendent vertices in a tree~$T$.

\begin{figure}[ht]
  \label{fig-delta}
  \center
  \includegraphics [width = 8cm]{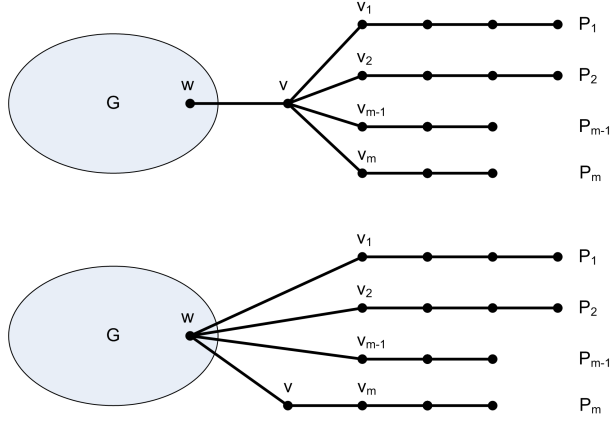}
  \caption { \textit{ $\delta$-transformation on vertex $v$. } }
\end{figure}

\begin{lemma} \label{1}
Let $T$ be a tree rooted at the center vertex $u$ with at least two
vertices of degree~3. Let $v \in \{z| \ deg(z) \geq 3, z \neq u\}$ be
a vertex with the largest distance $d (u, v)$ from the center vertex. Then for the
$\delta$-transformation tree $T' = \delta (T, v)$, it holds
$$
H (T') > H (T).
$$
\end{lemma}

\begin{proof}
We follow the symbols in Definition \ref{de-delta}.  Let $G$ be the
component of $T-wv$  containing the vertex $w$ (as shown in Fig. 3).
Let $R=\{P_1, P_2, \ldots, P_{m-1}\}$. After
$\delta$-transformation, the distances between vertices from $G$ and
$R$ decreased by $1$, while the distances between vertices from $R$
and $Q = P_m \cup \{v\}$ increased by $1$. By direct calculation, we
have
$$
H (T') - H (T) = \hspace{-0.3cm} \sum_{x \in G, y \in R} \frac{1}{d (x, y) - 1} - \frac{1}{d (x,
y)} \  + \sum_{x \in Q, y \in R} \frac{1}{d (x, y) + 1} - \frac{1}{d (x, y)}.
$$

According to the assumption, there is an induced path $P = w w_1 w_2
\ldots w_k$ in $G$, with length at least $\max \{n_1, n_2, \ldots,
n_m\}$. For each path $P_i$, $1 \leq i \leq m - 1$, it follows
\begin{eqnarray*}
D_i &=& \hspace{-0.3cm} \sum_{x \in G, y \in P_i} \left ( \frac{1}{d (x, y) - 1} - \frac{1}{d (x,
y)} \right )+
\sum_{x \in Q, y \in P_i} \left ( \frac{1}{d (x, y) + 1} - \frac{1}{d (x, y)}\right ) \\
&>& \hspace{-0.3cm} \sum_{x \in P, y \in P_i} \left  ( \frac{1}{d (x, y) - 1} - \frac{1}{d (x, y)}
\right ) -
\sum_{x \in Q, y \in P_i} \left (\frac{1}{d (x, y)} - \frac{1}{d (x, y) + 1} \right)\\
&\geq& 0.
\end{eqnarray*}

Hence
$$
H (T') - H (T) = \sum_{i = 1}^{m - 1} D_i > 0.
$$
Since $T$ contains at least two vertices of degree at least 3, we
have strict inequality. Therefore, if we move pendent paths $P_i$ towards the center vertex $u$ of $T$
along the path $P$, the Harary index increases.
\end{proof}

\begin{lemma}\label{2}\cite{gutman}
Let $G$ be a connected graph and $v \in V(G)$. Suppose that
$P=v_{0}v_{1}\ldots v_{k}, Q=u_{0}u_{1}\ldots u_{m}$ are two paths
with lengths $k$ and $m$ $(k\geq m\geq1)$, respectively. Let
$G_{k,m}$ be the graph obtained from $G, P, Q$ by identifying $v$
with $v_{0}$, $u$ with $u_{0}$, respectively. Then
$$
H(G_{k,m})>H(G_{k+1,m-1}).
$$
\end{lemma}
%

By Lemma \ref{1} and Lemma \ref{2}, we can get the following

\begin{proposition}\label{4}
Let $G_{0}$ be a connected graph and $u\in V(G_{0})$. Assume that
$G_{1}$ is the graph obtained from $G_0$ by attaching a tree $T$ ($T
\not \cong P_{k}$ and $T \not \cong S_{k}$) of order $k$ to $u$;
$G_{2}$ is the graph obtained from $G_0$ by identifying $u$ with an endvertex of a path $P_{k}$; $G_{3}$ is the graph obtained from $G_0$ by
 identifying $u$ with the center of a star $S_{k}$. Then
$$H(G_{2})<H(G_{1})< H(G_{3}).$$
\end{proposition}

Let $x = (x_1, x_2, \ldots, x_n)$ and $y = (y_1, y_2, \ldots, y_n)$
be two integer arrays of length~$n$. We say that $x$ majorizes $y$
and write $x \succ y$ if the elements of these arrays satisfy
following conditions:
\begin{enumerate}
\renewcommand{\labelenumi}{(\roman{enumi})}

\item $x_1 \geqslant x_2 \geqslant \ldots \geqslant x_n$ and $y_1 \geqslant y_2 \geqslant \ldots \geqslant
y_n$,
\item $x_1 + x_2 + \ldots + x_k \geqslant y_1 + y_2 + \ldots + y_k$,
for every $1 \leqslant k < n$,
\item $x_1 + x_2 + \ldots + x_n = y_1 + y_2 + \ldots + y_n$.
\end{enumerate}

\begin{theorem}
Let $p=(p_1, p_2, \ldots, p_k)$ and $q=(q_1, q_1, \ldots, q_k)$ be
two arrays of length $k \geqslant 2$, such that $p \prec q$ and $n =
p_1 + p_2 + \ldots + p_k = q_1 + q_2 + \ldots q_k$. Then
\begin{equation}
\label{eq:starlike} H (S (p_1, p_2, \ldots, p_k)) \geq H (S (q_1,
q_2, \ldots, q_k)).
\end{equation}
\end{theorem}

\begin{proof}
We will proceed by induction on the size of the array $k$. For $k =
2$, we can directly apply transformation from Lemma~\ref{2} on tree
$S (q_1, q_2)$ several times, in order to get $S (p_1, p_2)$. Assume
that the inequality (\ref{eq:starlike}) holds for all lengths less
than or equal to $k$. If there exist an index $1 \leqslant m < k$ such
that $p_1 + p_2 + \ldots + p_m = q_1 + q_2 + \ldots + q_m$, we can
apply the induction hypothesis  on two parts $S (q_1, q_2, \ldots,
q_m) \cup S (q_{m + 1}, q_{m + 2}, \ldots, q_k)$ and get $S (p_1, p_2,
\ldots, p_m) \cup S (p_{m + 1}, p_{m + 2}, \ldots, p_k)$.

Otherwise, we have strict inequalities $p_1 + p_2 + \ldots + p_m < q_1 + q_2 + \ldots
+ q_m$ for all indices $1 \leqslant m < k$. We can transform tree $S (q_1, q_2, \ldots, q_k)$ into
$$
S (q_1, q_2, \ldots, q_{s - 1}, q_{s} - 1, q_{s + 1}, \ldots, q_{r - 1}, q_{r} + 1, q_{r + 1}, \ldots, q_k),
$$
where $s$ is the largest index such that $q_1 = q_2 = \ldots = q_s$ and $r$ is the smallest index
such that $q_r = q_{r + 1} = \ldots = q_k$. The condition $p \prec q$ is preserved, and we can continue
until the array $q$ transforms into $p$, while at every step we increase the Harary index.
\end{proof}

\begin{corollary}\label{cor:order}
Let $T = S (n_1, n_2, \ldots, n_k) $ be a starlike tree with $n$
vertices and $k$ pendent paths. Then
$$
H (B_{n, k}) \leq H (T) \leq H (BS_{n, k}) .
$$
The left equality holds if and only if $T \cong B_{n,k}$ and the
right equality holds if and only if $T \cong BS_{n,k}$.
\end{corollary}

\section{Main results}
\vskip -0.1cm

\subsection{Trees with given number of pendent vertices}

Let ${\cal{T}}_{n,k}$ $(2\leq k\leq n-1)$ be the set of trees on $n$ vertices with $k$ pendent
vertices. If $k=2$, then the tree is just the path $P_n$. If $k=n-1$, then the tree is just the
star $S_n$. Therefore, we can assume that $3 \leq k \leq n-2$ in the sequel. It was proved in
\cite{IlIl09} that among $n$-vertex trees with given number $k$ of pendent vertices or given number
$q$ of vertices of degree two, $BS_{n, k}$ and $BS_{n, n - q - q}$ have minimal Wiener index.

\begin{theorem}\label{5}
Of all the trees on $n$ vertices  with $k$ $(3\leq k\leq n-2)$
pendent vertices, $BS_{n,k}$ is the unique tree having maximal
Harary index.
\end{theorem}

\begin{proof}
Suppose $T\in {\cal{T}}_{n,k}$ has the maximal Harary index, rooted
at the center vertex. Let $S_{T} = \{ v\in V(T) : deg(v)\geq 3 \}$.

If $|S_{T}| = 1$, then by Corollary \ref{cor:order}, it follows that
 $BS (n, k)$ is the unique tree that maximizes
Harary index. If $|S_{T}|\geq 2$, then there must exist at least two vertices of degree at
least 3 and there are only pendent paths attached below them.
We can consider $T$ as the rooted tree at the center vertex, and
choose the vertex $v$ of degree at least 3 that is furthest from the center vertex.
After applying $\delta$-transformation, we increase $H (G)$
while keeping the number of pendent vertices fixed -- which is
contradiction.
\end{proof}

\begin{theorem}\label{55}
Among trees with fixed number $q$  $(0\leq q\leq n-1)$ of vertices
of degree two, $T_{n,n-1-q}$ is the unique tree having maximal
Harary index.
\end{theorem}

\begin{proof}
The proof is similar to that of Theorem~\ref{5}.
%

Consider $\delta$-transformed tree $T' = \delta (T, v)$.
The vertex $v$ has degree greater than two, while the vertex $w$ has degree greater than or equal to two.
Among vertices on the pendent paths $P_1, P_2, \ldots, P_m$, there are
$$
S = (n_1 - 1) + (n_2 - 1) + \ldots + (n_m - 1) = \sum_{i = 1}^m n_i - m
$$
vertices of degree two.

If $w$ has degree two in $T$, then after $\delta$-transformation,
$v$ will have degree two, and the number of vertices of degree two in
$T'$ remains the same as in $T$. Otherwise, assume that $deg (w) > 2$. We can apply one
transformation from Lemma \ref{2}, and get new tree $T''$ with $m +
1$ pendent paths attached at vertex $w$ with lengths $n_1, n_2,
\ldots, n_m, 1$. This way we increased Harary index, while the
number of vertices with degree two in trees $T$ and $T''$ are the
same.

By repetitive application of these transformations, we prove that the starlike tree has maximal value
of Harary index among trees with $q$ vertices of degree two. The number of pendent paths is exactly
$k = \sum_{i = 1}^m n_i - q = n - 1 - q$, and by Corollary \ref{cor:order} it follows that
balanced starlike tree $BS (n, n - 1- q)$ is the unique tree that maximizes Harary index.
\end{proof}

By Lemma \ref{2}, we have the following chain of inequalities
\begin{equation}
\label{eq:order-balanced} H(P_{n}) = H(BS_{n,2}) < H(BS_{n,3}) < \ldots < H(BS_{n,n-1}) = H
(S_{n}).
\end{equation}

Notice that Lemma 2.2 from \cite{das}, $H (P_n) \leq \frac{(n + 2)(n - 1)}{4} = H (S_n)$ follows directly from these inequalities.

\subsection{Trees with given matching or independence number}

Du and Zhou in \cite{DuZh09} proved that among $n$-vertex trees with given matching number $m$,
the spur $A_{n, m}$ minimizes the Wiener index.

\begin{lemma}\label{6}
For arbitrary $\frac{n-1}{2} \leq m \leq n-1$, there holds
$$
H(A_{n,m})=\frac{1}{24} \left( 3 n^2+2 m n+m^2-9 m+19 n -22 \right).
$$
\end{lemma}

\begin{proof}
There are four types of vertices in the tree $A_{n, m}$. Denote with $D' (v)$ the sum of all
reverse distances from $v$ to all other vertices.
\begin{itemize}
\item For the center vertex, $D' (v) = \frac{m}{1} + \frac{n - m -
1}{2}$;
\item For each pendant vertex attached to the center vertex, $D' (v) = \frac{1}{1} + \frac{m - 1}{2} + \frac{n - m -
1}{3}$;
\item For each vertex of degree $2$, different from the center vertex, $D' (v) = \frac{2}{1} + \frac{m - 1}{2} + \frac{n - m -
2}{3}$;
\item For each pendant vertex not attached to the center vertex,
$D' (v) = \frac{1}{1} + \frac{1}{2} + \frac{m - 1}{3} + \frac{n - m - 2}{4}$.
\end{itemize}

After summing above contributions to the Harary index, we get
\begin{eqnarray*}
H(A_{n,m}) &=&  1 \cdot \frac{n + m - 1}{2} + (2m - n + 1) \cdot \frac{2n + m + 1}{6} \\
&& + \ (n - m - 1) \cdot \frac{2n + m + 5}{6} + (n - m - 1) \cdot \frac{3n + m + 8}{12} \\
&=& \frac{1}{24} \left( 3 n^2+m^2-22-9 m+19 n+2 m n \right).
\end{eqnarray*}
This proves the results.
\end{proof}

\begin{theorem}\label{7}
Let $T$ be a tree on $n$ vertices with matching number $\beta$. Then
$$
H(T)\leq \frac{1}{24} \left( 6 n^2 - 4 \beta n + \beta^2 + 9 \beta +
10 n  -22 \right),
$$
with equality holding if and only if $T\cong
A_{n,n-\beta}$.
\end{theorem}
\begin{proof}
Suppose that $T$ has  $k$  pendent vertices, then
$$k\leq \beta+n-2\beta=n-\beta,$$
and by Theorem \ref{5}, we have $H(T)\leq H(BS_{n,k})$. Using
equation (\ref{eq:order-balanced}), it follows $H(BS_{n,k}) \leq
H(BS_{n,n-\beta}) = H(A_{n,n-\beta})$ since $n-\beta \geq
\frac{n}{2}$. Finally, $H(T)\leq H(A_{n,n-\beta})$, with equality
holding if and only if $T\cong A_{n,n-\beta}$. By Lemma \ref{6}, we
obtain the explicit relation for $H(A_{n,n-\beta})$.
\end{proof}

By Theorem \ref{7}, we have the following corollary.
\begin{corollary} \label{8}
Let $T$ be a tree of order $n$ with perfect matching. Then
$$
H(T) \leq \frac{1}{4}(17n^2 + 58n -88),
$$
with equality holding if and only if $T\cong A_{n,\frac{n}{2}}$.
\end{corollary}

\begin{theorem}\label{9}
Let $T$ be a tree on $n$ vertices with independence number $\alpha$. Then
$$
H(T)\leq \frac{1}{24} \left( 3 n^2+2 \alpha n+\alpha^2-9 \alpha+19
n-22 \right),
$$
with equality holding if and only if $T\cong A_{n,\alpha}$.
\end{theorem}

\begin{proof}
Since all pendent vertices form an independent set, it follows $k\leq \alpha$. Every tree is
bipartite graph, and we get $\alpha \geq \lceil\frac{n}{2}\rceil$. By Theorem \ref{5}, we have
$H(T)\leq H(BS_{n,k})$ and $H(BS_{n,k})\leq H(BS_{n,\alpha})=H(A_{n,\alpha})$. Therefore, $H(T)\leq
H(A_{n,\alpha})$, with equality holding if and only if $T\cong A_{n,\alpha}$.
\end{proof}

%

\subsection{Trees with given diameter or radius}

Let $C_{n,d}(p_1, p_2, \ldots, p_{d-1})$ be a caterpillar on $n$ vertices obtained from a path
$P_{d+1}=v_0v_1\ldots v_{d-1}v_d$ by attaching $p_i \geq 0$ pendant vertices to $v_i$, $1\leq i
\leq d-1$, where $n=d+1+\sum_{i=1}^{d-1}p_i$. Denote
$$
C_{n,d,i}=C_{n,d}(\underbrace{0, \ldots, 0}_{i-1}, n-d-1, 0 \ldots ,0).
$$
Obviously, $C_{n,d,i}=C_{n,d, n-i}$. In \cite{IlIlSt09} and \cite{LiPa08} it is shown that
caterpillar $C_{n,d, \lfloor d / 2\rfloor}$ has maximal Wiener index among trees with fixed
diameter $d$.


\begin{theorem}
Among trees on $n$ vertices and diameter $d$,
$C_{n,d, \lfloor d / 2\rfloor}$ is the unique tree having maximal Harary index.
\end{theorem}

\begin{proof}
Let $T$ be $n$-vertex tree with diameter $d$ having maximal Harary index. Let $P_{d+1}=v_0v_1\ldots
v_d$ be a path of length $d$. By Proposition~\ref{4}, all trees attached to the path $P_{d+1}$ must
be stars, which implies that $T \cong C_{n,d}(p_1, p_2, \ldots , p_{d-1})$. Applying Lemma~\ref{1}
at those vertices of degree at least 3 in $C_{n,d}(p_1, p_2, \ldots , p_{d-1})$, we get that $T
\cong C_{n,d,i}$ for some $1 \leq i \leq d - 1$. By Lemma \ref{2} we get that the maximal Harary
index is achieved uniquely for $i=\lfloor\frac n2\rfloor$. This completes the proof.
\end{proof}

For a tree $T$ with radius $r(T)$, it holds $2 r (T) = d (T)$ or $2r (T) - 1 = d (T)$. Using
transformation from Lemma \ref{2} applied to a center vertex, it follows that $H (C_{n, 2r,\lfloor
d/ 2 \rfloor}) < H (C_{n, 2r-1,\lfloor  d/ 2 \rfloor})$.

\begin{corollary}
Let $T$ be a tree on $n$ vertices with radius $r \geq 2$. Then
$$
H (T) \leq H (C_{n, 2r-1,\lfloor  d/ 2 \rfloor}),
$$
with equality if and only if $T \cong C_{n, 2r-1,\lfloor  d/ 2 \rfloor}$.
\end{corollary}

If $d>2$, we can apply the transformation from Lemma~\ref{2} at the
center vertex in $C_{n,d, \lfloor d/2\rfloor}$  to obtain
$C_{n,d-1, \lfloor (d-1)/2 \rfloor}$. Thus,
$$
H(P_{n}) = H(C_{n,n-1,\lfloor (n-1)/2 \rfloor}) < \ldots < H(C_{n,3,1})< H(C_{n,2,1}) = H (S_{n}).
$$

Also, it follows that $H(C_{n,3,1})$ has the second maximum Harary
index among trees on $n$ vertices.

\subsection{Trees with given maximum vertex degree}

Chemical trees (trees with maximum vertex degree at most four) provide the graph representations of
alkanes \cite{GuPo86}. It is therefore a natural problem to study trees with bounded maximum
degree. It was proven in \cite{DoEn01} that among $n$-vertex trees with the maximum degree
$\Delta$, the broom $B_{n, \Delta}$ has maximal Wiener index.

\begin{lemma}\label{77}
For arbitrary $2 \leq \Delta \leq n - 1$, there holds
$$
H(B_{n,\Delta}) = n \cdot H_{n - \Delta} - n + \Delta +
\frac{(\Delta - 1)(\Delta - 2)}{4} + \frac{\Delta - 1}{n - \Delta +
1},
$$ where $H_k=1+\frac 12+\ldots+ \frac 1k $ be the $k$-th harmonic
number.
\end{lemma}

\begin{proof}
The Harary index of $B_{n,\Delta}$ can be calculated as the sum of
Harary index of $P_{n - \Delta + 1}$, the sum of reverse distances
between $\Delta - 1$ pendent vertices and vertices from a long path,
and the sum of reverse distance between pendent vertices.
\begin{eqnarray*}
H(B_{n,\Delta}) &=& H (P_{n - \Delta + 1}) + \frac{ \binom{\Delta - 1}{2} }{2} + (\Delta - 1) H_{n - \Delta + 1} \\
&=& n \cdot H_{n - \Delta} - n + \Delta + \frac{(\Delta - 1)(\Delta
- 2)}{4} + \frac{\Delta - 1}{n - \Delta + 1}.
\end{eqnarray*}
\end{proof}

\begin{theorem}
Let $T $ be a tree on $n$ vertices with the maximum degree~$\Delta$.
Then $H (T)\leq H (B_{n,\Delta})$. The equality holds if and only if
$T \cong B_{n,\Delta}$.
\end{theorem}

\begin{proof}
Fix a vertex $v$ of degree $\Delta$ as a root and let $T_1, T_2, \ldots, T_{\Delta}$ be the trees
attached at~$v$. By Proposition \ref{4}, all subtrees attached to $u_i$ are paths for $1 \leq i\leq
\Delta$, the Harary index increases. This implies the result.
\end{proof}

If $\Delta>2$, we can apply the transformation from Lemma~\ref{2} at the
vertex of degree~$\Delta$ in $B_{n, \Delta}$ and obtain
$B_{n, \Delta-1}$. Thus,
$$
H(S_{n}) = H(B_{n,n-1}) > H(B_{n,n-2}) > \ldots > H(B_{n,3}) >
H(B_{n,2})= H (P_{n}).
$$

Also, it follows that $B_{n, 3}$ has the second minimum Harary index
among trees on $n$ vertices.

\section{Concluding remarks}
\vskip -0.1cm

In this paper, we presented the partial ordering of starlike trees based on the Harary index and we
derived the trees with the second maximal and the second minimal Harary index. We characterized the
extremal trees with maximal Harary index and fixed number of pendent vertices, the number of
vertices of degree two, matching number, independence number, radius and diameter. In addition, we
characterized the extremal trees with minimal Harary index and given maximum degree. We concluded
that in the all presented classes, the trees with maximum values of Harary index are exactly those
trees with the minimal Wiener index $W (G)$, and vice versa.

The complete $\Delta$-ary tree is defined as follows. Start with the root having $\Delta$~children.
Every vertex different from the root, which is not in one of the last two levels, has exactly
$\Delta - 1$~children. In the last level, while not all nodes have to exist, the nodes that do
exist fill the level consecutively. Thus, at most one vertex on the level second to last has its
degree different from $\Delta$ and~$1$.


In~\cite{GuFMG07} the authors proposed these trees to be called
\emph{Volkmann trees}, as they represent alkanes with minimal Wiener
index \cite{FiHo02}. The computer search among trees with up to 24
vertices reveals that the complete $\Delta$-ary trees attain the
maximum values of $H (G)$ among the trees with the maximum vertex
degree~$\Delta$.

\begin{conjecture}
For any $k \geqslant 2$, the complete $\Delta$-ary tree has maximum
value of $H (G)$ among trees on $n$ vertices with maximum
degree~$\Delta$.
\end{conjecture}

It would be interesting for further research to consider the extremal unicyclic and bicyclic graphs
with respect to Harary index, and compare the results with those for the Wiener index.

\end{document}